\documentclass[1p]{elsarticle}

\usepackage{lineno, hyperref}
\usepackage{amsmath}
\usepackage{amssymb}
\usepackage{amsfonts}
\usepackage{graphicx}
\usepackage{graphics}
\usepackage{epstopdf}
\usepackage{caption}
\usepackage{comment}
\usepackage{csquotes}
\usepackage{float}
\usepackage{enumerate}
\usepackage{latexsym}
\usepackage{subcaption}
\usepackage[toc]{appendix}
\usepackage{xcolor}
\usepackage{cleveref}

\journal{}

\newtheorem{remark}{Remark}
\biboptions{sort&compress}
\begin{document}

\begin{frontmatter}

\title{Approximations for modeling light scattering by spheres with uncertainty in physical parameters}

\author[]{Akif Khan}
\ead{khanakif91@gmail.com}

\author[]{Murugesan Venkatapathi}
\ead{murugesh@iisc.ac.in}

\address{Computational  and  Statistical  Physics  Laboratory,  Department  of  Computational  and  Data Sciences, Indian Institute of Science, Bangalore - 560012}

\date{}

\begin{abstract}
Uncertainty in physical parameters can make the solution of forward or inverse light scattering problems in astrophysical, biological, and atmospheric sensing applications, cost prohibitive for real-time applications. For example, given a probability density in the parametric space of dimensions, refractive index and wavelength, the number of required evaluations for the expected scattering increases dramatically. In the case of dielectric and weakly absorbing spherical particles (both homogeneous and layered), we begin with a Fraunhofer approximation of the scattering coefficients consisting of Riccati-Bessel functions, and reduce it into simpler nested trigonometric approximations. They provide further computational advantages when parameterized on lines of constant optical path lengths. This can reduce the cost of evaluations by large factors $\approx$ 50, without a loss of accuracy in the integrals of these scattering coefficients. We analyze the errors of the proposed approximation, and present numerical results for a set of forward problems as a demonstration.
\end{abstract}


\end{frontmatter}

\newcommand{\D}{\mathrm{d}}
\section{Introduction:}
Light scattering as an indicator of the optical properties of particles has been used widely for sensing and diagnosis in biological, atmospheric, astrophysical, and combustion applications. In biological applications, it is frequently used in measurement of optical properties of tissues \cite{tuchin1997light}, bacterial colonies \cite{Bae2019colonies}, and medical diagnostics by automated cell classification or sorting \cite{Doornbos96trapcell, Venkatapathi2008scattering}. Similarly, there are many atmospheric and astrophysical applications that involve an estimation of particle sizes and refractive indices \cite{jones1988modelling, min2005modeling}. These particles could be aerosols in the atmosphere or cosmic grain dust in the interstellar medium. With modern computational tools, modeling light scattering from a particle with known physical parameters has become relatively trivial, even when they are irregularly shaped \cite{Mazeron96scattering, dunn1997finite, Maltsev2011scattering}.

Evaluation using known analytical solutions is significantly more efficient compared to solving the electromagnetic scattering problem using numerical models, and such analytical solutions are known only for (homogeneous and layered) spheres, infinite cylinders and spheroids \cite{Latimer1975ellipsoids, bohren2008absorption, Arun2016Sommerfeld}. Consequently, numerical approaches for solution such as T-matrix \cite{waterman1965matrix, mackowski1996scattering}, Discrete Dipole Approximation(DDA) \cite{draine1994discrete} and the Finite Difference Time Domain(FDTD) \cite{dunn1997finite} have gained prominence over the years in addressing other geometries, notwithstanding their higher cost of computation.

In the case of uncertainty in the physical parameters, one may estimate the expected scattering from a large number of particles representing a distribution of physical properties, and this is called a \textit{forward} problem. When these particles are irregular in shapes and/or random in orientations, they may be replaced by a statistical ensemble of particles like spheres and spheroids for which the closed form solution for the scattering are known \cite{Asano1980RandomSpheroids, jones1988modelling, min2005modeling}. On the other hand, \textit{inverse} problems estimate the distribution of the properties of the particles, given a finite number of scattering measurements \cite{Markel2018BRT, Markel2019inverse}. A reasonably accurate solution to inverse problems, may also require the solution of such forward problems multiple times, increasing the computational cost further \cite{Burger2001levelset, Cakoni2013InverseIntegral, Kern2016Inverse}. Thus, fast algorithms to evaluate the scattering for layered and homogeneous spheres having a distribution in size and index becomes even more imperative. Well-known approximations in light scattering from a single particle include the RGD (Rayleigh-Gans-Debye), WKBJ (Wentzel–Kramers–Brillouin-Jeffreys) and the Born approximations, and they are based on physical insights into the problem \cite{Nikolay1999ScatApprox, Bereza2017Born}. The Born approximation has been found to be useful in reducing computations for inverse problems in medical imaging and radar based identification, for example \cite{Trattner2007BornImaging, Li2010BornRadars}.

For particles of a fixed shape and orientation with respect to the incident light, having a probability distribution in size and refractive index, the expected intensity of light can be evaluated using the integral:
\begin{equation}
    I=\iint_S p(\hat{x},\hat{m})\sigma_{sca}(\hat{x},\hat{m}) \D^l \hat{x} \D^l\hat{m} \label{eqn:eqn 1}
\end{equation}
where $\hat{x}=[x_1,\dots,x_l],\hat{m}=[m_1,\dots,m_l]$. $p$ and $\sigma_{sca}$ are the probability density function and scattering cross-sections respectively, evaluated in the parametric space of size parameters $x$ (i.e. size normalized by wavelength of interest) and refractive indices $m$ of a non-homogeneous particle. $\sigma_{sca}$ may represent the cross-sections of extinction, total scattering, the back-scattering, the angular differential scattering, or the forward scattering. For homogeneous ($l=1$) and layered non-homogeneous ($l > 1$) spherical particles, they are evaluated using the scattering coefficients $a_n$ and $b_n$ derived from the boundary conditions of the scattering problem for each orthogonal mode numbered $n$. Thus a viable evaluation of the above integrals depends on (a) efficient computation of the integrand given by the scattering coefficients, and (b) reducing the number of evaluations of the integrand, required for the numerical integration. The latter is done using suitable numerical methods for integration, and is not the subject of this work.

Here, it is highlighted that any periodicity of the scattering cross-sections can be exploited in simplifying its evaluation, even when this periodicity in the parametric space is non-trivial. A less known \textit{circular law} for dielectric spheres which constrains the scattering coefficients to a circle in the complex plane, is first introduced to motivate the oscillatory nature of the scattering coefficients with varying physical parameters. Later, the Fraunhofer approximations of the Riccati-Bessel functions are used to reduce the scattering cross-sections into trigonometric expressions. These approximations of the scattering coefficients should be distinguished from the first and second order expansions of the Green functions for evaluating the diffraction by an obstacle, known as Fraunhofer and Fresnel diffractions respectively. Using nested forms of these trigonometric expressions and parameterizing the evaluations on lines of constant optical lengths, the number of trigonometric functions to be evaluated are further reduced to a few. Since the magnitude of errors in the proposed approximation of the integrand is small, and the cancellation of resulting error due to its oscillatory signs, one can significantly reduce computational effort for the integrals such as \cref{eqn:eqn 1} without a loss of accuracy.

\section{A circular law of scattering coefficients for homogeneous and layered spheres:}\label{sec:circular law}

The extinction cross-section $C_{ext}$ of a particle, typically interpreted as the area of its \textit{shadow}, describes the reduction of intensity of light in the forward direction due to the absorption and scattering by the particle i.e. $C_{ext}=C_{sca}+C_{abs}$. For homogeneous and layered spheres, they are evaluated using \cite{bohren2008absorption}:
\begin{align}
    C_{sca}&=\frac{2\pi}{k^2}\sum_{n=1}^\infty (2n+1)(|a_n|^2 + |b_n|^2) \label{eqn:eqn 3}\\
    C_{ext}&=\frac{2\pi}{k^2}\sum_{n=1}^\infty (2n+1)(Re(a_n + b_n)) \label{eqn:eqn 4}
\end{align}
where $a_n$ and $b_n$ are the scattering coefficients for each mode numbered by integers $n$, and they are functions of the size parameter $x$=$kR$ and the relative refractive index $m$ with respect to the background medium. Here, $R$ is the radius of the sphere and the wave number $k$=$2\pi/\lambda$, where $\lambda$ is the wavelength of the incident plane wave. For homogeneous spheres, the scattering coefficients are given in the form of Riccati-Bessel functions $\psi_n$ and $\chi_n$ as \cite{bohren2008absorption}: 
\begin{align}
    a_n &= \frac{m\psi_n(mx)\psi_n'(x)-\psi_n(x)\psi_n'(mx)}{m\psi_n(mx)\xi_n'(x)-\xi_n(x)\psi_n'(mx)} \label{eqn:a_homogeneous}\\
    b_n &= \frac{\psi_n(mx)\psi_n'(x)-m\psi_n(x)\psi_n'(mx)}{\psi_n(mx)\xi_n'(x)-m\xi_n(x)\psi_n'(mx)}\label{eqn:b_homogeneous}
\end{align}
where $\xi_n = \psi_n + i\chi_n$. For a sphere with an additional layer with the outer size parameter $y$=$kR_2$ and refractive index $m_2$, in addition to an inner sphere given by $x$ and $m_1$, the corresponding scattering coefficients are:

\begin{align}
    a_n& = \frac{\psi_n(y) \left[ \psi'_n\left(m_2y\right) - A_n\chi'_n\left(m_2y\right)\right] - m_2\psi'_n(y) \left[ \psi_n\left(m_2y\right) - A_n\chi_n\left(m_2y\right)\right]}{\xi_n(y) \left[ \psi'_n\left(m_2y\right) - A_n\chi'_n\left(m_2y\right)\right] - m_2\xi'_n(y) \left[ \psi_n\left(m_2y\right) - A_n\chi_n\left(m_2y\right)\right]}\label{eqn:eqn 7}\\
    b_n&= \frac{m_2\psi_n(y) \left[ \psi'_n\left(m_2y\right) - B_n\chi'_n\left(m_2y\right)\right] - \psi'_n(y) \left[ \psi_n\left(m_2y\right) - B_n\chi_n\left(m_2y\right)\right]}{m_2\xi_n(y) \left[ \psi'_n\left(m_2y\right) - B_n\chi'_n\left(m_2y\right)\right] - \xi'_n(y) \left[ \psi_n\left(m_2y\right) - B_n\chi_n\left(m_2y\right)\right]}\label{eqn:eqn 8}
\end{align}
where,
\begin{align}
    A_n&=\frac{m_2 \psi_n\left(m_2x\right)\psi'_n\left(m_1x\right) - m_1 \psi'_n\left(m_2x\right)\psi_n\left(m_1x\right)}{m_2 \chi_n\left(m_2x\right)\psi'_n\left(m_1x\right) - m_1 \chi'_n\left(m_2x\right)\psi_n\left(m_1x\right)}\label{eqn: eqn 9}\\
    B_n&=\frac{m_2 \psi_n\left(m_1x\right)\psi'_n\left(m_2x\right) - m_1 \psi_n\left(m_2x\right)\psi'_n\left(m_1x\right)}{m_2 \chi'_n\left(m_2x\right)\psi\left(m_1x\right) - m_1 \psi'_n\left(m_1x\right)\chi_n\left(m_2x\right)}\label{eqn: eqn 10}
\end{align}
For dielectric spheres, the absorption cross-section $C_{abs}=0$ and we have $C_{ext}=C_{abs}$. Equating their expressions: 
\begin{equation}
    \sum_{n=1}^\infty (2n+1)(|a_n|^2 + |b_n|^2)=\sum_{n=1}^\infty (2n+1)(Re(a_n + b_n))
    \label{eqn:eqn 11}
\end{equation}
Now, the above equality is true for any sphere of arbitrary size parameters, refractive indices and number of layers in the sphere. This generality allows us to assume that the $a_n$'s and $b_n$'s are independent variables in equating the individual terms of the series on both the sides \cite{garg2017fano}. Thus we get the ansatz:
\begin{equation}
    Re(a_n)=|a_n|^2 \text{ and } Re(b_n)=|b_n|^2 \label{eqn:eqn 12}
\end{equation}
Let both $a_n$ and $b_n$ be complex numbers of the form $z=x+iy$, then the equation $Re(z)=|z|^2$ describes a circle in the complex plane having radius $1$ and centered at $\left(\frac{1}{2},0\right)$, as:
\begin{align*}
    x^2+y^2-x=0 \implies
    \left(x-\frac{1}{2}\right)^2 + y^2 = \frac{1}{4}
    \label{eqn:eqn 14}
\end{align*}
We call this ansatz as the \textit{circular law} of scattering coefficients of layered and homogeneous dielectric spheres. A more rigorous proof of this law for homogeneous spheres was presented by van de Hulst  many decades ago \cite{hulst1981light}.
\begin{figure}[H]
  \centering
  \includegraphics[width=1.0\linewidth]{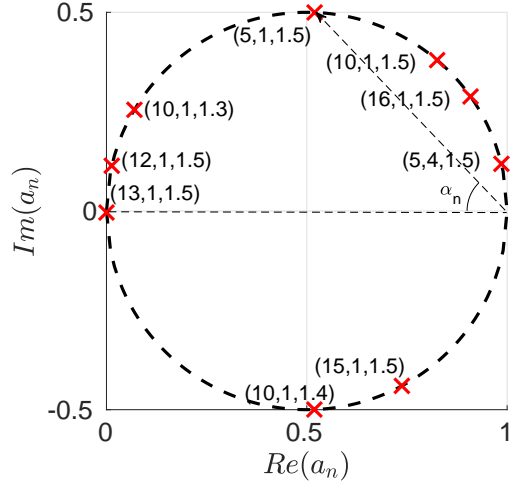}
  \caption{The circular law for dielectric spheres is demonstrated for the scattering coefficient $a_n$ of homogeneous spheres, with random values of triplets ($x$, $n$, $m$). The law holds true for an arbitrary size parameter $x$, mode number $n$, and refractive index $m$, for any number of layers $l$ in the particle.}
  \label{fig:figure 1}
\end{figure}

\section{Pseudo-periodicity of scattering coefficients and their trigonometric approximations:}
Since the scattering coefficients of spheres have been restricted to a circle in the complex plane, it is convenient to represent them by a single variable given by an angle. We shall now define the angles $\alpha$ and $\beta$, where $\alpha$ is $\frac{\pi}{2}-\text{Arg}(a_n)$ and $\beta$ is $\frac{\pi}{2}-\text{Arg}(b_n)$. Conveniently, $\alpha_n$ and $\beta_n$ are given by the relations, $|\sin\alpha_n| = |a_n|$ and $|\sin \beta_n| = |b_n|$. Thus the scattering cross-section becomes:
\begin{equation}
    C_{sca}(m,x)=\frac{2\pi}{k^2}\sum_{n=1}^\infty (2n+1)(\sin^2\alpha_n + \sin^2\beta_n)
    \label{eqn:eqn 15}
\end{equation}

It is known that the Riccati-Bessel functions are also given by \cite{abramowitz1999ia}:
\begin{align}
\psi(\rho)&=P(n+\frac{1}{2},\rho)\sin(\rho-\frac{1}{2}n\pi)+Q(n+\frac{1}{2},\rho)\cos(\rho-\frac{1}{2}n\pi)\label{eqn:eqn 20}\\
\chi(\rho)&=(-1)^{n+1}\left(P(n+\frac{1}{2},\rho)\cos(\rho+\frac{1}{2}n\pi)-Q(n+\frac{1}{2},\rho)\sin(\rho+\frac{1}{2}n\pi)\right)\label{eqn:eqn 21} 
\end{align}
where the terms P, Q are given by the sums of a finite series:
\begin{align}
\begin{split}
    P(n+\frac{1}{2},\rho)&=1-\frac{(n+2)!}{2!\Gamma(n-1)}(2\rho)^{-2}+\frac{(n+4)!}{4!\Gamma(n-3)}(2\rho)^{-4}-...\\
    &=\sum_{k=0}^{\lfloor\frac{n}{2}\rfloor}(-1)^k(n+\frac{1}{2},2k)(2\rho)^{-2k}\label{eqn:eqn 22}
\end{split}\\
\begin{split}  
    Q(n+\frac{1}{2},\rho)&=\frac{(n+1)!}{1!\Gamma(n)}(2\rho)^{-1}-\frac{(n+3)!}{3!\Gamma(n-2)}(2\rho)^{-3}+\frac{(n+5)!}{5!\Gamma(n-4)}(2\rho)^{-5}-...\\
    &=\sum_{k=0}^{\lfloor\frac{n-1}{2}\rfloor}(-1)^k(n+\frac{1}{2},2k+1)(2\rho)^{-2k-1} \label{eqn:eqn 23}
\end{split}
\end{align}

\begin{align*}
(n+\frac{1}{2},k)&=\frac{(n+k)!}{k!\Gamma(n-k+1)} \mbox{\hspace{1.5in}(n=0,1,2...)}
\end{align*}
The asymptotic forms of the Riccati-Bessel functions for large arguments(for $\rho\gg n^2$), called the Fraunhofer approximations, are given by just the first term,
\begin{align}\label{eqn:psi_defn}
    \psi_n(\rho)&=\sin(\rho-n\frac{\pi}{2})\\\label{eqn:chi_defn}
    \chi_n(\rho)&=(-1)^{n+1}\cos(\rho+n\frac{\pi}{2})
\end{align}
and they can be used for exact evaluations of the scattering coefficients only for smaller mode numbers $n$ and relatively larger spheres ($x \gg n^2$).

\subsection{Approximation for homogeneous spheres}
Using $\xi_n = \psi_n + i\chi_n$, we first rewrite \cref{eqn:a_homogeneous} for $a_n$ in terms of a numerator $N$ consisting only of $\psi_n$ and its derivative, and a $D$ containing $\chi_n$ as well. Note that $N$ and $D$ are real numbers for a real refractive index $m$.
\begin{equation}\label{eqn:N_D}
a_n = \frac{N}{N+iD} \text{ and } |a_n|^2=\sin^2\alpha_n=\frac{\tan^2\alpha_n}{1+\tan^2\alpha_n} \implies \tan \alpha_n = -\frac{N}{D}
\end{equation}
Given the Fraunhofer approximations of $\psi(\rho)$ and $\chi(\rho)$ in \cref{eqn:psi_defn} and \cref{eqn:chi_defn}, one can reduce the evaluations for $\sin^2\alpha_n$ of homogeneous spheres to simpler nested trigonometric functions. We first demonstrate this using an expression for $\sin^2\alpha_1$. Applying Fraunhofer approximations for $\psi$ and $\chi$ in $N$ and $D$,
\begin{equation}\label{eqn:tan_alpha1}
    -\tan\alpha_1 \approx \frac{-mx\cos(mx)\sin(x) + x\cos(x)\sin(mx)}{mx\cos(mx)\cos(x)+x\sin(x)\sin(mx)}
\end{equation}

\begin{figure}[H]
  \begin{subfigure}{0.9\textwidth}
  \centering
  \includegraphics[width=0.8\linewidth]{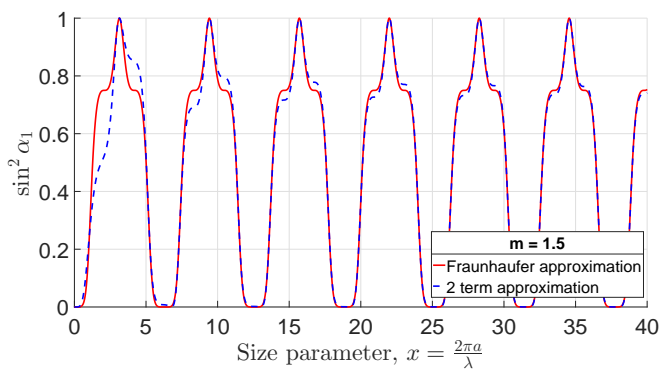} 
  \caption{Fraunhofer and two-term approximations with a constant refractive index $m$ and varying size parameter $x$. Note the sharp nature of the resonances.}
\end{subfigure}
\begin{subfigure}{0.9\textwidth}
  \centering
  \includegraphics[width=0.8\linewidth]{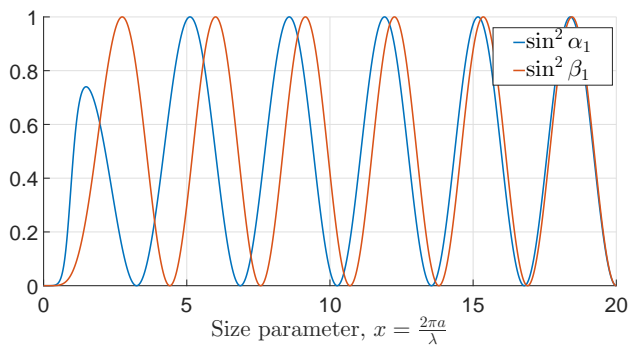}
  \caption{The exact Bessel function evaluations of $\sin^2\alpha_1$ and $\sin^2\beta_1$ for a given $mx=c=20$, illustrates the slowly varying $\delta$ in the approximation $\sin^2(x-\delta)$, converging to $\delta=x$ when $x=c$ and $m=1$. It begins with a small high index sphere on the left varying towards a large low index sphere on the right.}
\end{subfigure}
\caption{Homogeneous spheres}
\label{fig:figure 2}
\end{figure}

The above approximation reduces the evaluation of the scattering coefficients from Bessel functions to trigonometric evaluations, and their plots are compared with the two-term (Fresnel) approximations of the scattering coefficients in Figure \ref{fig:figure 2}a. We do not present the explicit expressions of the Fresnel approximations of the scattering coefficients here, for brevity, and it can be seen that the contribution of the higher order terms are small for larger arguments. From the approximations of $\sin^2\alpha_1$ varying with $x$ and a constant $m$, we observe that although it is roughly periodic it cannot be represented in the form of any elementary trigonometric function. Note that we are only interested in values of $x$ larger than the mode numbers i.e. $x \gtrsim n$, where such evaluations of scattering coefficients are required \cite{bohren2008absorption}. For particles with $x < 1$, the mode-less Rayleigh scattering approximations would suffice. Thus, our interest is in the Fresnel regime $n \leq x \leq cn^2$,\cite{Heitman2015Bessel} where $c$ is $\mathcal{O}(1)$, and the Fraunhofer approximation of scattering coefficients are shown to be reasonable for integrals, and it can be reduced further to a simple trigonometric expression for faster evaluations. Dividing \cref{eqn:tan_alpha1} by $mx\cos(mx)\cos(x)$:
\begin{equation*}
    \tan\alpha_1 \approx \frac{\tan(x)-\frac{\tan(mx)}{m}}{1+\frac{\tan(x)\tan(mx)}{m}}=\tan(x - \delta)
\end{equation*}
and thus,
\begin{equation}\label{eq:homogeneous_approx_an}
    \sin^2\alpha_1 \approx \sin^2(x - \delta)
\end{equation}
where
\begin{equation}\label{eqn:final_homogeneous}
    \delta = \tan^{-1}(\frac{\tan(mx)}{m})
\end{equation}
Similarly,
\begin{equation}\label{eq:homogeneous_approx_bn}
    \sin^2\beta_1 \approx \sin^2(x - \gamma)
\end{equation}
where
\begin{equation}
    \gamma = \tan^{-1}(m\tan(mx))
\end{equation}
We can further reduce the approximation in \cref{eq:homogeneous_approx_an} to two trigonometric evaluations using constant optical path lengths $mx = c$. Note that when the refractive index $m$=1, the particle should vanish, and it is satisfied by these approximations. For larger values of $n$, the approximations in \cref{eq:homogeneous_approx_an} and \cref{eq:homogeneous_approx_bn} repeat themselves with an odd and even alternation between $\sin^2 \alpha_n$ and $\sin^2 \beta_n$. This is due to the odd and even behavior exhibited by the Fraunhofer approximation of the Riccati-Bessel functions, given in equations \eqref{eqn:psi_defn} and \eqref{eqn:chi_defn}.

\begin{remark}
For any given mode `n' of the homogeneous sphere, the Fraunhofer approximation of $|a_n|^2$ and $|b_n|^2$ reduces to $\sin^2(x-\delta)$ where $\delta$ is a trigonometric function of size parameter `x' and refractive index `m'.
\end{remark}

The exact Bessel forms of $\sin^2\alpha_1$ and $\sin^2\beta_1$ in Figure \ref{fig:figure 2}b, varying with $x$ but a constant $mx$, show the approximate periodicity indicated above by the trigonometric approximations. This periodicity is exploited for reducing the computational cost in evaluating the scattering due to a single particle. A more detailed description of the errors in this Fraunhofer approximation is presented in section \ref{sec:Errors}. This apparent loss of accuracy by applying the Fraunhofer approximation in the Fresnel regime comes at the gain of a significant reduction in computing. The number of trigonometric evaluations in the approximation is reduced to \textit{two} compared to the \textit{eight} Bessel functions in the full evaluations. As computation of a spherical Bessel function costs more than a trigonometric function, 
we expect a significant decrease in the computational cost in the evaluation of the scattering cross-sections in the parametric space.

\subsection{Approximation for layered spheres}

\begin{figure}[h]
  \begin{subfigure}{0.9\textwidth}
  \centering
  \includegraphics[width=0.8\linewidth]{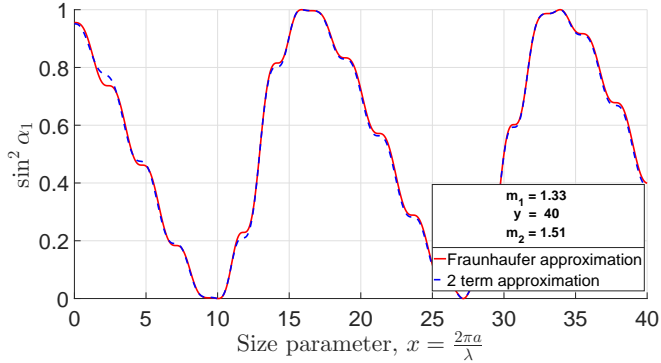} 
  \caption{Fraunhofer and two-term approximations of a layered sphere for a varying inner radius $a$, with a constant refractive index and outer radius of the sphere. The variation is periodic but the pattern is non-trivial.}
\end{subfigure}
\begin{subfigure}{0.9\textwidth}
  \centering
  \includegraphics[width=0.8\linewidth]{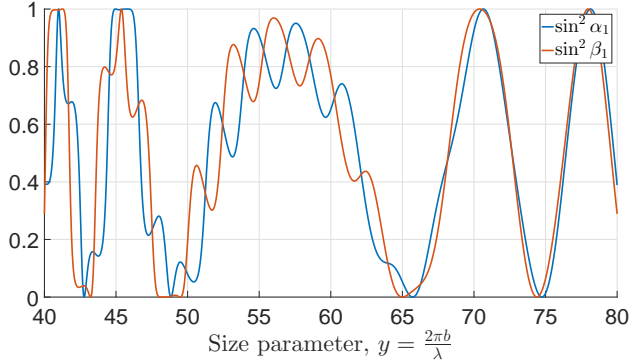} 
  \caption{$\sin^2\alpha$ and $\sin^2\beta$ are non-periodic for a constant $m_2y$ and a varying outer size parameter $y$. The distinct patterns changing from left to right represents the variation of the outer sphere from a thin high index layer to a relatively large low index layer. Here $m_1$=1.33 and $x$=40.}
\end{subfigure}
\caption{Layered spheres}
\label{fig:figure 3}
\end{figure}

For layered spheres, just as in the case of homogeneous spheres, varying one of the variables such as the inner radius $R_1$ ($x=kR_1$) with constant parameters $m_1$, $m_2$ and $y=kR_2$ does not present a simple periodic behavior of the scattering constants, as shown in Figure \ref{fig:figure 3}a. Even a constant optical size of the outer sphere $m_2y=c$, with a varying $y$, does not present viable trigonometric approximations, as shown in Figure \ref{fig:figure 3}b. This can also be inferred using the more cumbersome expression for $\tan\alpha_n$ of the layered spheres. Using $\xi_n = \psi_n + i\chi_n$, we first rewrite \cref{eqn:eqn 7} for $a_n$ in terms of a numerator $N$ consisting only of $\psi_n$ and its derivative, and a $D$ containing $\chi_n$ as well. Note that $N$ and $D$ are real numbers for real refractive indices $m$.
\begin{equation*}
a_n = \frac{N}{N+iD} \text{ and } |a_n|^2=\sin^2\alpha_n=\frac{\tan^2\alpha_n}{1+\tan^2\alpha_n} \implies \tan \alpha_n = -\frac{N}{D}
\end{equation*}
\begin{equation}
    -\tan\alpha_n =  \frac{\psi_n(y) \left[ \psi'_n\left(m_2y\right) - A_n\chi'_n\left(m_2y\right)\right] - m_2\psi'_n(y) \left[ \psi_n\left(m_2y\right) - A_n\chi_n\left(m_2y\right)\right]}{\chi_n(y) \left[ \psi'_n\left(m_2y\right) - A_n\chi'_n\left(m_2y\right)\right] - m_2\chi'_n(y) \left[ \psi_n\left(m_2y\right) - A_n\chi_n\left(m_2y\right)\right]}  
\end{equation}

Define $P$, $Q$ as 
\begin{align*}
    P &= \psi'_n\left(m_2y\right) - A_n\chi'_n\left(m_2y\right)\\
    Q &= \psi_n\left(m_2y\right) - A_n\chi_n\left(m_2y\right)
\end{align*}

Now applying the Fraunhofer approximations of \cref{eqn:psi_defn} and \cref{eqn:chi_defn}, $\tan\alpha_1$ becomes:
\begin{equation*}
    \tan\alpha_1 \approx \frac{y\cos y P + m_2y\sin y Q}{-y\sin y P + m_2y\cos y Q}
\end{equation*}
Dividing the numerator and the denominator by $m_2yQ\cos y$:
\begin{equation*}
    \tan\alpha_1 \approx \frac{\frac{P}{m_2Q} + \tan y}{1 - \frac{P\tan y}{m_2Q}}=\tan(y + \delta_1)
\end{equation*}
where $\delta_1$ is given by
\begin{equation*}
    \tan\delta_1 = \frac{P}{m_2Q}
\end{equation*}
Thus we have,
\begin{equation*}
    \alpha_1 \approx y + \delta_1 = y + \tan^{-1}\left(\frac{P}{m_2Q}\right)
\end{equation*}
Also $\frac{P}{Q}$ in the Fraunhofer approximation becomes
\begin{equation*}
    \frac{P}{Q} \approx \frac{\sin(m_2y)+A_1\cos(m_2y)}{-\cos(m_2y)+A_1\sin(m_2y)}
\end{equation*}

Dividing by $\cos (m_2y)$, we get
\begin{equation*}
    -\frac{P}{Q} \approx \frac{\tan(m_2y) + A_1}{1 - A_1\tan(m_2y)}
\end{equation*}
Thus $-\frac{P}{Q} \approx \tan(m_2y - \delta_1')$
where $\tan\delta_1' = -A_1$. Recall that to find $a_n$ for layered spheres, we need to evaluate $A_n$ which is given by:
\begin{equation*}
    A_n = \frac{m_2\psi_n(m_2x)\psi_n'(m_1x)-m_1\psi_n'(m_2x)\psi_n(m_1x)}{m_2\chi_n(m_2x)\psi_n'(m_1x)-m_1\chi_n'(m_2x)\psi_n(m_1x)}
\end{equation*}
We make the transformations $\frac{m_1}{m_2}=M$, $m_2x=X$ we get $A_n$ as 
\begin{equation*}
    A_n = \frac{M\psi_n(MX)\psi_n'(X)-\psi_n(X)\psi_n'(MX)}{M\psi_n(MX)\chi_n'(X)-\chi_n(X)\psi_n'(MX)}
\end{equation*}
which is just the expression for $-\tan\alpha_n$ for homogeneous spheres of size parameter $X$ and index $M$. Thus, recalling its trigonometric approximation in \cref{eqn:final_homogeneous}:
\begin{equation*}
    \delta_1' \approx X-\tan^{-1}(\frac{\tan MX}{M}) = m_2x - \tan^{-1}\left(\frac{m_2}{m_1}\tan (m_1x)\right)
\end{equation*}

Thus $\alpha_1$ becomes:
\begin{equation}\label{eqn:final_layered}
    \alpha_1 \approx y - \tan^{-1}\left(\frac{\tan (m_2y - \delta_1')}{m_2}\right)
\end{equation}

Using the above approximations of $\alpha_1$ we get: 
\begin{equation*}
    \sin^2\alpha_1 \approx \sin^2\left(y - \tan^{-1}\left(\frac{\tan (m_2y - \delta_1')}{m_2}\right)\right)
\end{equation*}

By similar approximations,
\begin{equation*}
    \sin^2\beta_1 \approx \sin^2\left(y - \tan^{-1}\left(\tan(m_2[m_2y-\gamma_1'])\right)\right)
\end{equation*}
with $\gamma_1' = m_1x - \tan^{-1}\left(\frac{m1}{m2}\tan(m_1x)\right)$. For larger mode numbers $n > 1$, the above approximations for $\sin^2\alpha_1$ and $\sin^2\beta_1$ repeat themselves with an odd and even alternation between them. This is due to the odd and even behavior exhibited by the Fraunhofer approximation of the Bessel function, given in equations \eqref{eqn:psi_defn} and \eqref{eqn:chi_defn}. Thus we have reduced the scattering cross-section of the layered sphere in the form of Bessel functions, to approximate trigonometric expressions, for the ease of computing.

\begin{remark}
For any given mode `n' of the layered sphere, the Fraunhofer approximation of $|a_n|^2$ and $|b_n|^2$ reduces to $\sin^2(y-\delta)$ where $\delta$ is a nested trigonometric function of inner size parameter `x' and refractive index `$m_1$', outer size parameter `y' and refractive index `$m_2$'.
\end{remark}

\subsection{Errors in the approximation} \label{sec:Errors}
Before we plot the error and speedup obtained by the trigonometric approximations for practical forward problems involving integrals, we analyze the errors of the trigonometric approximations in evaluation of the integrand. This allows us to see that the errors are oscillatory in sign and thus involve significant cancellations. 

\begin{figure}
  \centering
  \includegraphics[width=0.8\linewidth]{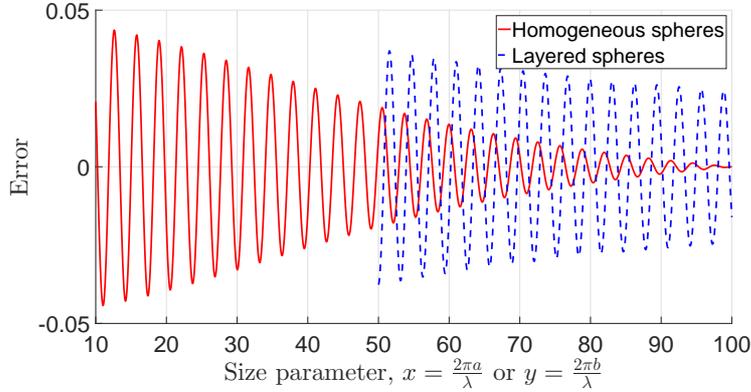}
  \caption{The error in trigonometric approximation is bounded and oscillatory, and thus the integral error (in the next figure) can be expected to approach zero for sufficiently large intervals. In the above plot for $n=1$, the parameters are for homogeneous spheres: $c = mx = 100$, for layered spheres: $c_3 = m_2y= 100$, $c_2 = m_2x = 30$, $c_1 = m_1x = 50$.}
  \label{fig:figure 4}
\end{figure}

Recall that we have rewritten the scattering coefficients as $|a_n|=|\sin\alpha|$ where the angle $\alpha$ represents the angle $\frac{\pi}{2}-\text{Arg}(a_n)$, and $a_n$ lies on a circle in complex plane. The same applies to spherical coefficients $b_n$. The corresponding error in the integral of $|a_n|^2$ over a parametric space $S$ when using the Fraunhofer approximation $\alpha_f$ for $\alpha$ is:
\begin{equation*}
    \int_S \left(\sin^2\alpha - \sin^2\alpha_f\right)ds=\int_S \sin(\alpha+\alpha_f)\sin(\alpha-\alpha_f) ds
\end{equation*}
where we have used trigonometric identities to rewrite the integrand into a product of two components. As the Fraunhofer approximation $\alpha_f \longrightarrow \alpha$ for larger values of $x$, the error vanishes as expected. More significantly, the above trigonometric expression highlights the oscillatory nature of magnitude and sign of the error, as $\alpha$ spans the whole range when $a_n$ traverses the circle in the complex plane as described in section \ref{sec:circular law}. Thus, while the error tends to zero when $x \gg n^2$ (see Figure \ref{fig:figure 4}), the errors of this approximation are oscillatory for $n \leq x \leq n^2$.

\begin{figure}
  \centering
  \includegraphics[width=0.8\linewidth]{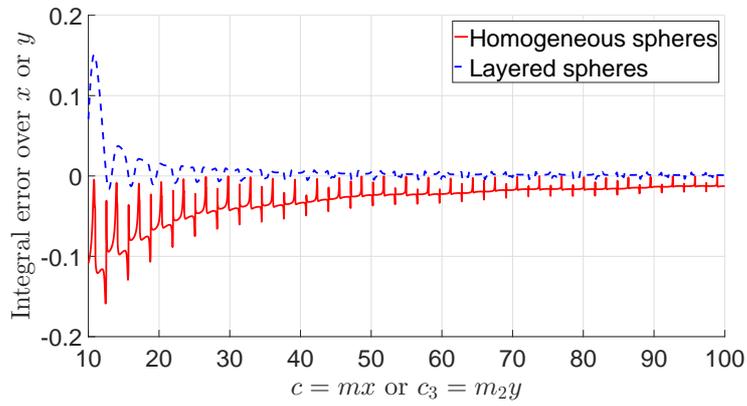}
  \caption{Due to the oscillatory nature of the error with respect to $x$, the error for $n=1$ in the integral from $0$ to $x$ is much lower.}
  \label{fig:figure 5}
\end{figure}

\begin{figure}
  \centering
  \includegraphics[width=0.8\linewidth]{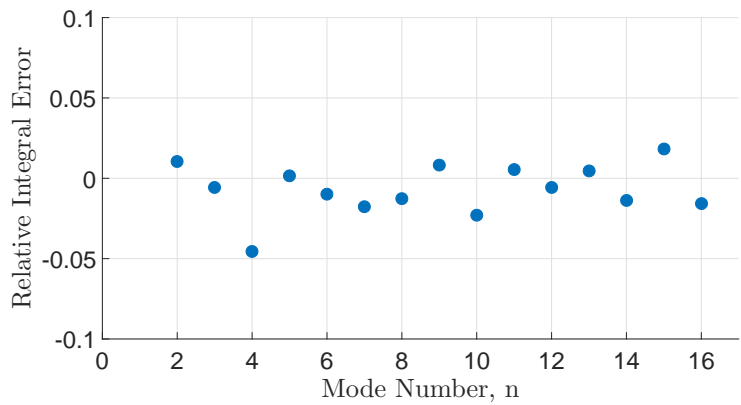}
  \caption{The relative integral error was computed for each mode $n$, over the domain $n \leq x\leq 2\pi n$ for optical path lengths $mx = 2\pi n$ of a homogeneous sphere.}
  \label{fig:figure 6}
\end{figure}

Due to the oscillatory nature of the error, we expect the errors to cancel out for a reasonably large domain of quadrature, and indeed we see that in Figure \ref{fig:figure 5} along a line in the parametric space with a constant optical length. Integrals over a two-dimensional parametric space involve stronger cancellations resulting in even smaller relative errors, as we demonstrate in the next section. Figure \ref{fig:figure 6} highlights that the relative errors in the integrals when using the Fraunhofer approximation in the Fresnel regime, are indeed low for any mode number $n$ in general.

\subsection{Extension to back-scattering and forward-scattering evaluations}
The differential back-scattering and the differential forward scattering cross-sections are given by\cite{garg2017fano}:

\begin{equation}
    \sigma_b=\left[\frac{1}{(2k)^2}\left|\sum_{n=1}^{\infty}(2n+1)(-1)^n(a_n-b_n)\right|^2\right]
    \label{eqn:eqn 28}
\end{equation}

\begin{equation}
    \sigma_f=\left[\frac{1}{(2k)^2}\left|\sum_{n=1}^{\infty}(2n+1)(a_n+b_n)\right|^2\right]
    \label{eqn:eqn 29}
\end{equation}
From the circular law, we know that $|a_n|=|\sin\alpha_n|$ and $|b_n|=|\sin\beta_n|$. The expressions for $a_n$ and $b_n$ in terms of $\sin\alpha_n$ and $\sin\beta_n$ are:
\begin{eqnarray}
    a_n &=&|\sin\alpha_n|(\sin\alpha_n + i\cos\alpha_n)
    \label{eqn:eqn 30}\\
    b_n &=&|\sin\beta_n|(\sin\beta_n + i\cos\beta_n)
    \label{eqn:eqn 31}
\end{eqnarray}
Using the reduced Fraunhofer approximations discussed so far, thus, we can approximate coefficients $a_n$ and $b_n$ for evaluating the forward and back-scattering cross-sections as well.

\section{Numerical examples:}
We compared the accuracy of the above  trigonometric approximation with respect to the exact Bessel function evaluation of $C_{sca}$ for curves in the parametric space given by a constant optical length. Now we extend this to the entire parametric space of interest given by a few probability distributions and show that the above approximations are indeed both efficient and accurate. We use Gaussian quadrature to compute the integral $I$ in \cref{eqn:eqn 1}, and at the quadrature points we use both the above trigonometric approximations (in \cref{eqn:final_homogeneous} and \cref{eqn:final_layered}) and the full Bessel function evaluations of $C_{sca}(m,x)$ (using equations \ref{eqn:a_homogeneous} - \ref{eqn: eqn 10}), for direct estimations of the speed-up. We also present the relative errors in the integrals when we use the trigonometric approximations for evaluating the integral $I$ using the quadrature. The true values of the integral are given by the asymptotic values for a large number of points in the numerical quadrature where the integrand is evaluated using the exact Bessel function forms of $a_n$ and $b_n$.

\subsection{Variations of speedup over parametric points}

	\begin{figure}
		\parbox{7cm}{
			\includegraphics[width = 7 cm]{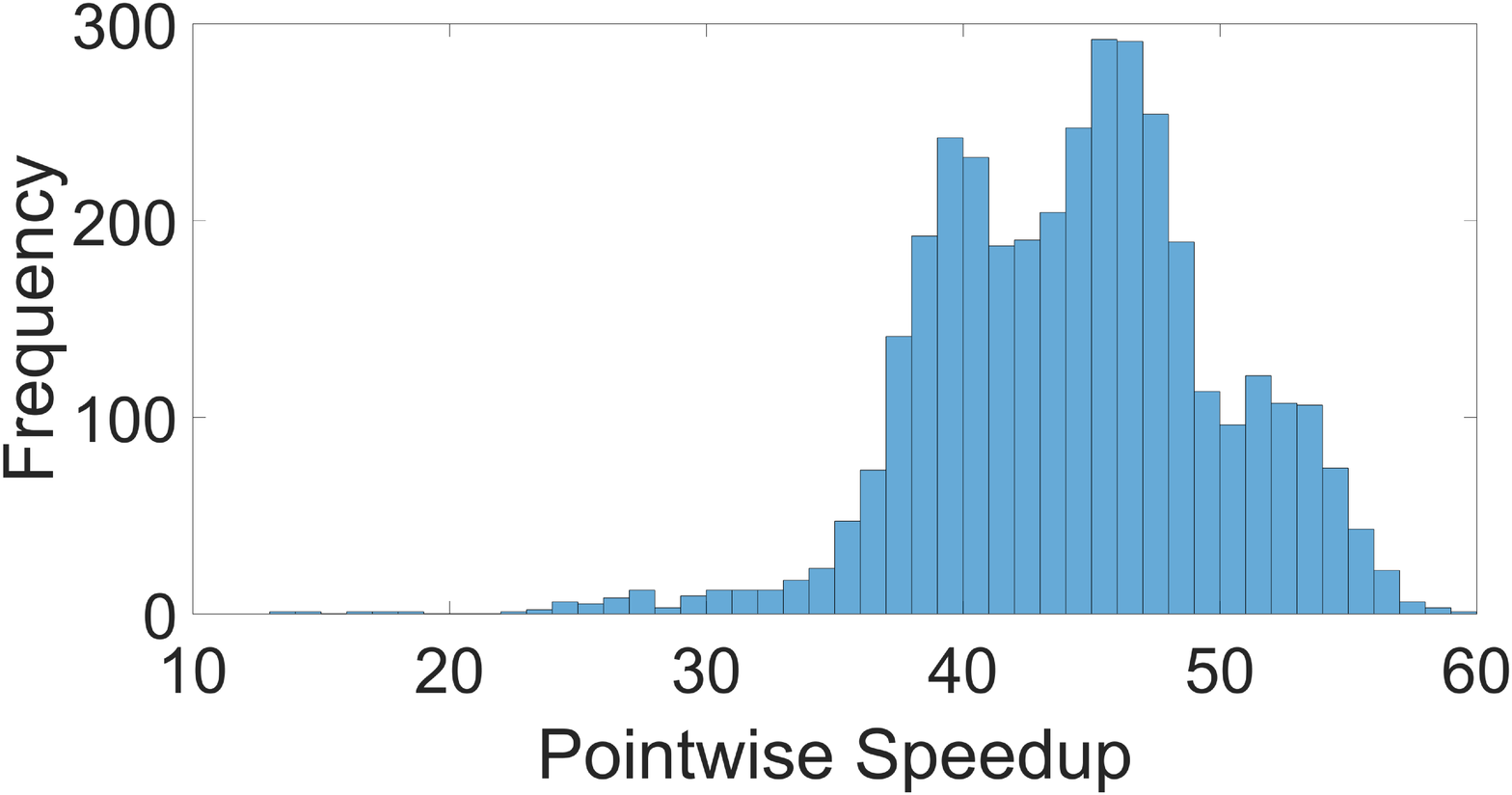}
		}
		\parbox{7cm}{
			\includegraphics[width = 7 cm]{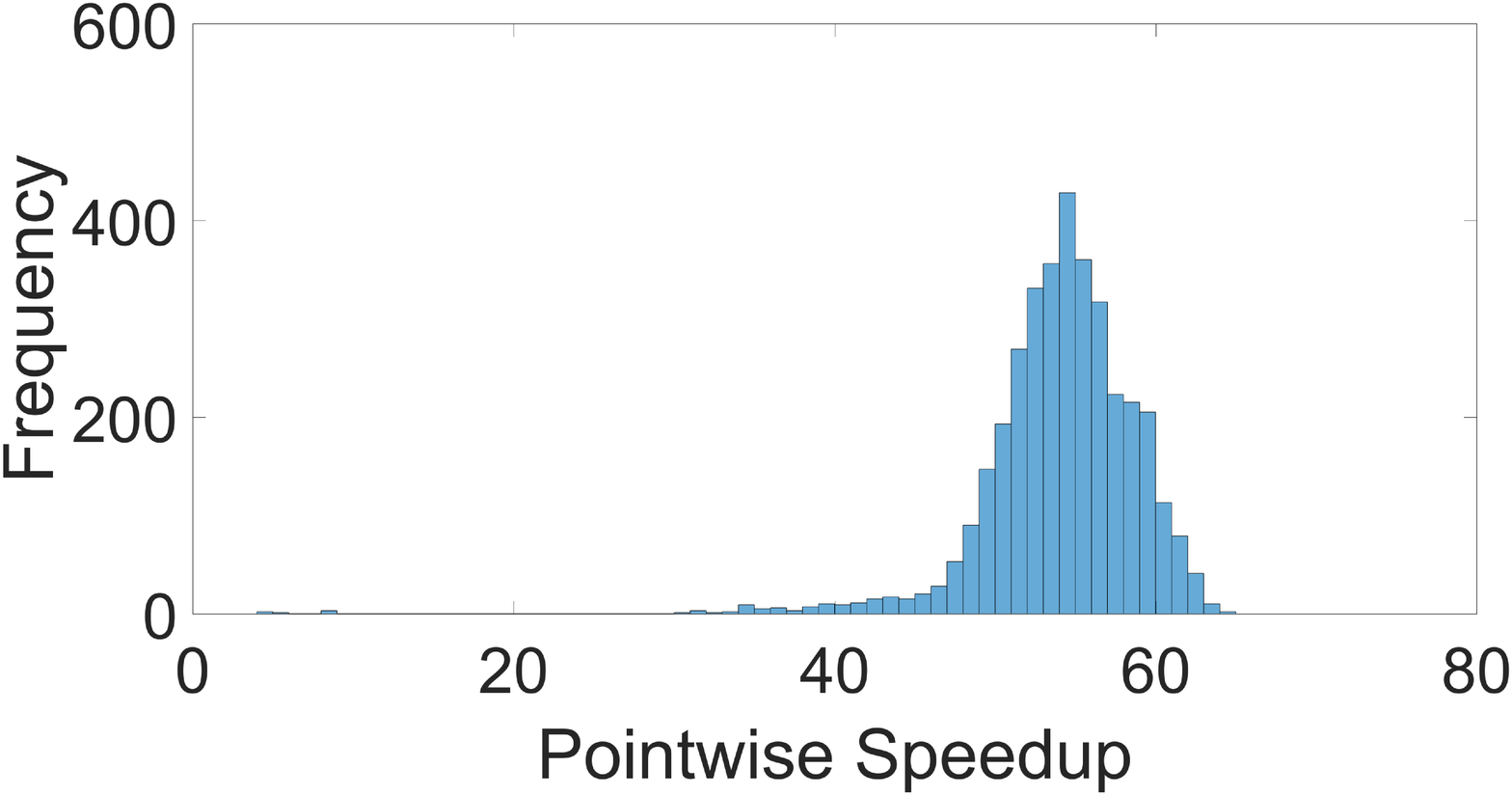}
		}
		\caption{The speedup for 3600 different points marginally varies with the parameters $m$ and $x$ (or $y$ for a layered sphere) as the convergence of the Bessel function in the exact evaluations depend upon the arguments. All function calls were invoked in MATLAB. \textit{Left}: Homogeneous spheres. \textit{Right}: Layered spheres.}\label{fig:figure 7}
	\end{figure}
We note that the speedup due to the trigonometric approximation of the scattering coefficients is not constant over the domain. It varies with $m$ and $x$ as the convergence of the Riccati-Bessel function in the exact evaluations depend upon its arguments. We define the point-wise speedup as the ratio of the wall-time taken to compute the scattering for a point in the parametric space i.e. a single particle with certain properties, using both the full evaluations and the trigonometric approximation.

The histograms in Figure \ref{fig:figure 7} show the frequency of the point wise speedups for homogeneous and layered spheres using 3600 ($60 \times 60$) points in $m$ and $x$. These results show that the cost of 8 Bessel functions for a homogeneous sphere and further arithmetic operations can cost more than 50 times the evaluation using the trigonometric approximation. Similarly, the trigonometric approximation for a layered sphere costs higher, but so does its full evaluation which consists of 16 Bessel functions and their arithmetic operations. This reduction in computing can be attained in principle for any number of layers of a radially inhomogeneous sphere each defined by a pair of values $m$, $x$.

\subsection{Uniform distribution}

\begin{figure}[H]
  \begin{subfigure}{0.9\textwidth}
  \centering
  \includegraphics[width=0.8\linewidth]{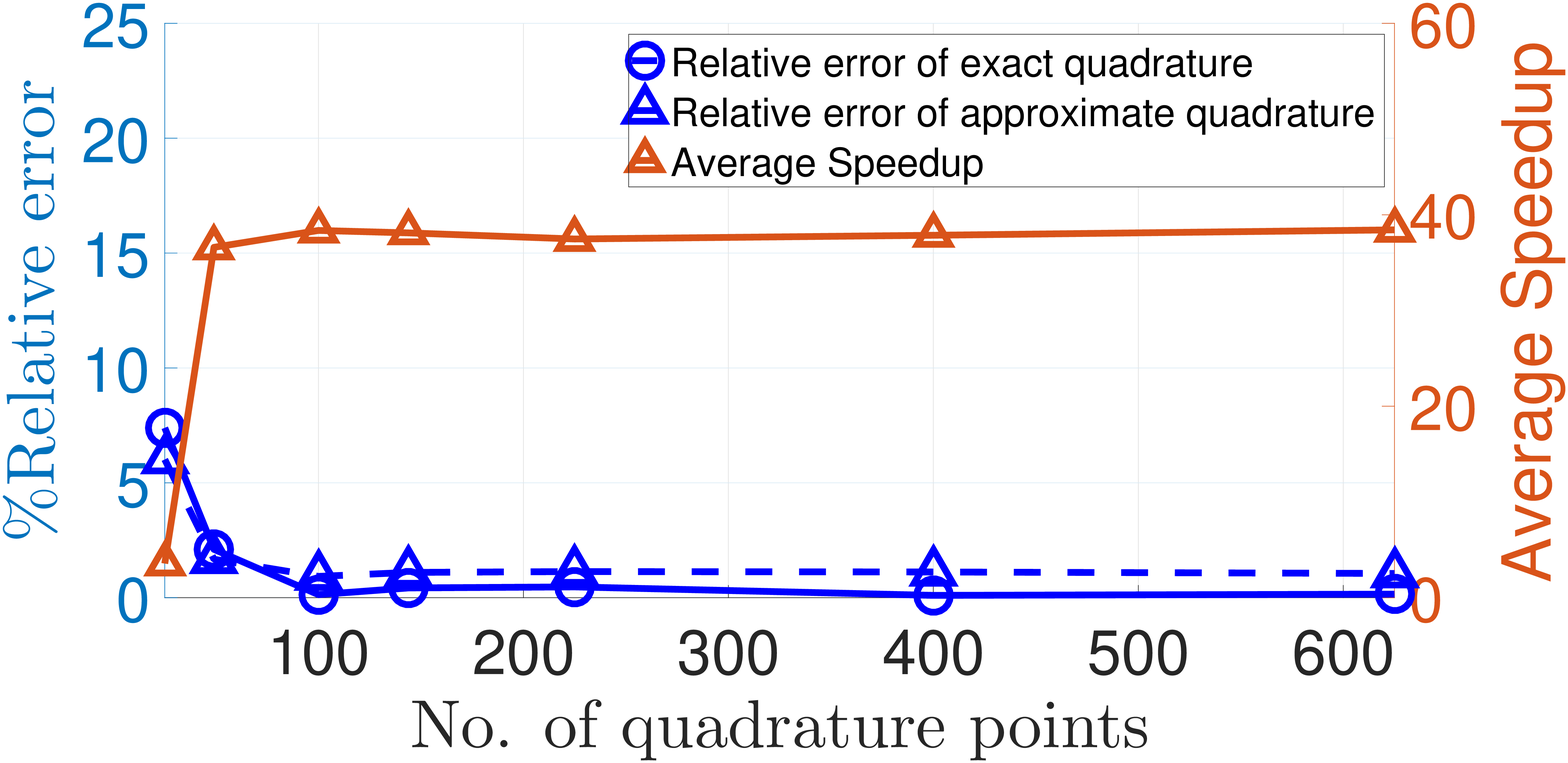}
  \caption{Homogeneous spheres}
\end{subfigure}
\begin{subfigure}{0.9\textwidth}
  \centering
  \includegraphics[width=0.8\linewidth]{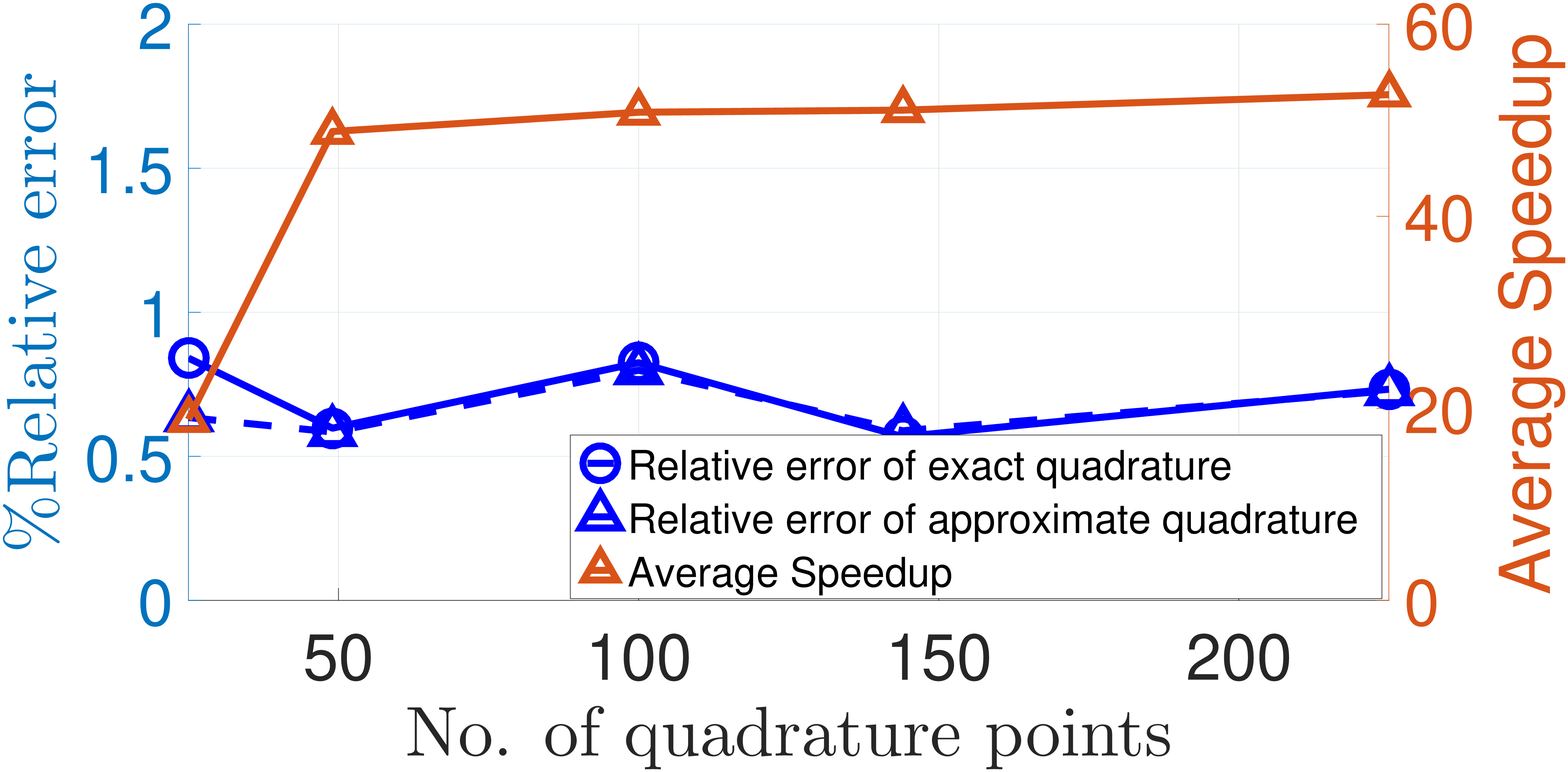} 
  \caption{Layered spheres}
\end{subfigure}
\caption{The error due to trigonometric approximation converges to that of the exact evaluations for a uniform distribution. The speedup so achieved is $\approx$ 60 for homogeneous as well as layered spheres.}
\label{fig:figure 8}
\end{figure}
The following values of the parameters $m$ and $x$ were considered and these ranges are motivated by the practical applications as discussed in in the introduction.
\begin{enumerate}
    \item Homogeneous spheres: $10 \leq x \leq 20$, $1.2 \leq m \leq 1.8$.\\
    Homogeneous spheres in this range are encountered in biological and combustion diagnostic applications. 
    \item Layered spheres: $40 \leq x \leq 60$, $1.25 \leq m_1 \leq 1.4$;
    We consider a PDF in $m_1$ and $x$ while having a constant thickness of outer shell as $y$=$x+20$ and refractive index $m_2$=1.51. Such problems represent atmospheric applications with soot deposited water droplets.
\end{enumerate}
Along with a description of the convergence, we also plot the average speed-up which reflects the comparative time required to compute the integrand using the modes $n=1,2,3$ for the original Bessel-function form and the trigonometric approximation of $a_n$ and $b_n$.

\subsection{Normal distribution}
We consider a two dimensional normal distribution for the size parameter $x$ and index $m$ of homogeneous and layered spheres to show the efficacy of the approximations for non-uniform distributions. With $\sigma$ and $\mu$ as the variance and mean respectively, the following normal distributions were simulated reflecting the range of parameters and applications described in the previous sub-section on uniform distribution.
\begin{enumerate}
    \item Homogeneous spheres: $\mu_x=15$, $\sigma_x=1.67$ and $\mu_m=1.5$, $\sigma_m=0.1$. 
    \item Layered spheres: $\mu_x=50$, $\sigma_x=3.33$ and $\mu_{m1}=1.325$, $\sigma_{m1}=0.025$. $y$=$x+20$ and refractive index $m_2$=1.51.
\end{enumerate}

\begin{figure}[H]
  \begin{subfigure}{0.9\textwidth}
  \centering
  \includegraphics[width=0.8\linewidth]{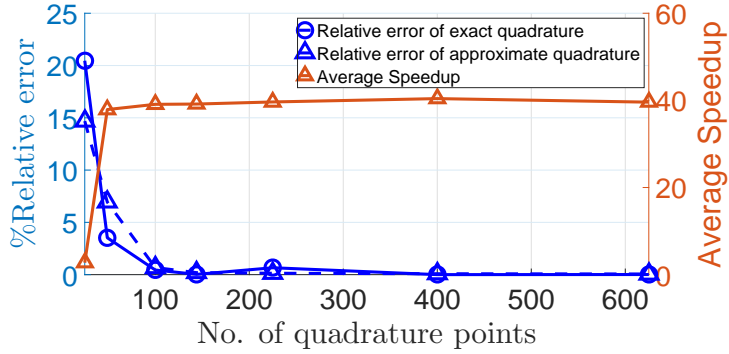}
  \caption{Homogeneous spheres}
\end{subfigure}
\begin{subfigure}{0.9\textwidth}
  \centering
  \includegraphics[width=0.8\linewidth]{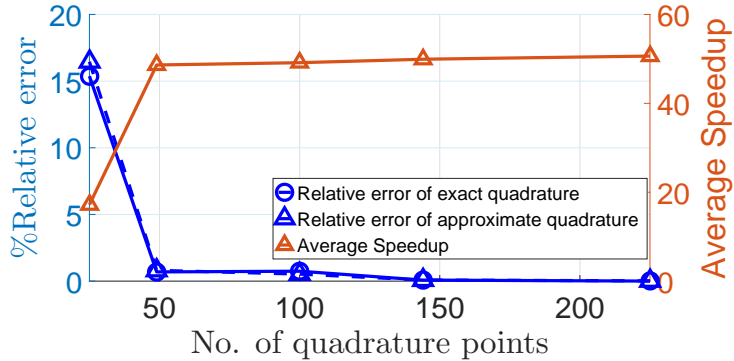} 
  \caption{Layered spheres}
\end{subfigure}
\caption{The error and speedup is seen to be distribution independent, and the error due to trigonometric approximation converges to that of the exact evaluations for a normal distribution.}
\label{fig:figure 9}
\end{figure}

\subsection{Bimodal distribution}
We can construct a bimodal distribution $p(x,m)$ for homogeneous or layered spheres using a weighted sum of two normal distributions.
\begin{equation}
    p(x,m)=w_1f_1(x,m) + w_2f_2(x,m)
\end{equation}
Where, $w_1$, $w_2$ are the weights such that $w_1 +w _2=1$ and $f_1(x,m)$, $f_2(x,m)$ are two normal distributions such that $|\mu_1-\mu_2|\geq 2\sigma$, where $\mu_1$, $\mu_2$ are the means of the two normal distributions and $\sigma$ is the average variance of the two normal distributions. This criteria is necessary to obtain two peaks (modes) in the bimodal distribution. Numerical experiments used equal weights $w_1$=$w_2$=0.5, for both homogeneous and layered spheres. For a demonstration we use the following values for $\sigma$ and $\mu$ reflecting applications referred in the previous sub-section on uniform distribution. 
\begin{enumerate}
    \item Homogeneous spheres: $f_1$ - $\mu_{x}=13$, $\sigma_{x}=1$ and $\mu_{m}=1.4$, $\sigma_{m}=0.06$.
    $f_2$ - $\mu_{x}=17$, $\sigma_{x}=1$ and $\mu_{m}=1.6$, $\sigma_{m}=0.06$.
    \item Layered spheres: $y$=$x+20$ and refractive index $m_2$=1.51. $f_1$ - $\mu_{x}=45$, $\sigma_{x}=3.33$ and $\mu_{m_1}=1.30$, $\sigma_{m_1}=0.02$.
    $f_2$ - $\mu_{x}=55$, $\sigma_{x}=3.33$ and $\mu_{m_1}=1.35$, $\sigma_{m_1}=0.02$.
\end{enumerate} 

\begin{figure}[H]
  \begin{subfigure}{0.9\textwidth}
  \centering
  \includegraphics[width=0.8\linewidth]{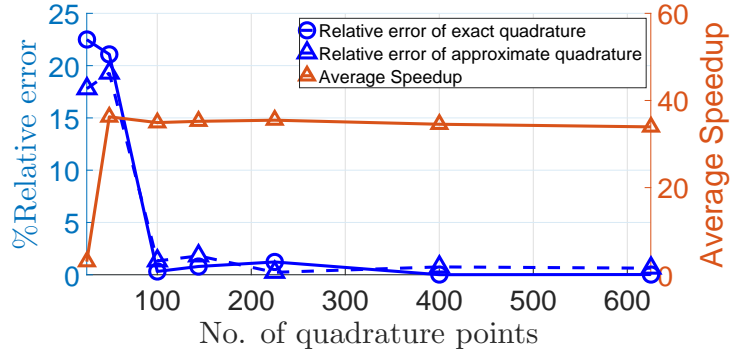}
  \caption{Homogeneous spheres}
\end{subfigure}
\begin{subfigure}{0.9\textwidth}
  \centering
  \includegraphics[width=0.8\linewidth]{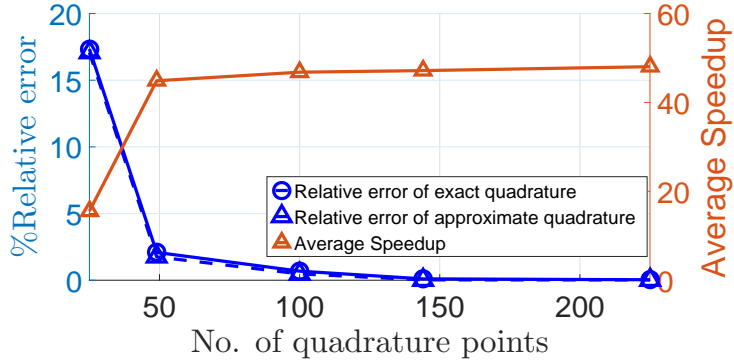} 
  \caption{Layered spheres}
\end{subfigure}
\caption{The error and speedup is shown to be distribution independent, and the error due to trigonometric approximation converges to that of the exact evaluations for even a bimodal distribution.}
\label{fig:figure 10}
\end{figure}

\bibliography{mybibfile}

\end{document}